\documentclass[10pt,a4paper,reqno]{amsart}
\usepackage{hyperref}

\usepackage[margin=1in,marginparwidth=0.8in]{geometry}
\usepackage{tikz}
\usepackage{graphicx}
\usepackage{todonotes}
\usepackage{mathtools}
\usepackage{shuffle}
\usepackage{color}
\usepackage{comment}

\usepackage[font=small]{caption}

\usepackage{ amssymb, mathrsfs }
\usepackage{bbm}
\usepackage{dsfont}

\def\={\:{}={}\:}

\DeclareMathOperator{\Li}{Li}

\newcommand{\Qb}{\mathbf{Q}}
\newcommand{\Mf}{\mathfrak{M}}

\newcommand{\IIr}[1]{I_{#1}^{N}}

\newtheorem{Thm}{Theorem}
\newtheorem{Cor}[Thm]{Corollary}
\newtheorem{Prop}[Thm]{Proposition}

\theoremstyle{definition}

\newtheorem{Rem}[Thm]{Remark}

\title[Functional equations of polygonal type for MPLs]{Functional equations of polygonal type for multiple polylogarithms in weights 5, 6 and 7}
\date{\today}
\author[Charlton]{Steven Charlton}
\address{Fachbereich Mathematik (AZ), Universit\"at Hamburg, Bundesstra\textup{\ss}e 55, 20146 Hamburg, Germany}
\email{steven.charlton@uni-hamburg.de}

\author[Gangl]{Herbert Gangl}
\address{Department of Mathematical Sciences, Durham University, Durham DH1 3LE, United Kingdom}
\email{herbert.gangl@durham.ac.uk}

\author[Radchenko]{Danylo Radchenko}
\address{ETH Zurich, Mathematics Department, Zurich 8092, Switzerland}
\email{danradchenko@gmail.com}


\begin{document}
	\begin{abstract}
	We present new functional equations in weights 5, 6 and 7 and use them for explicit depth reduction of multiple polylogarithms. These identities generalize the crucial identity~$\Qb_4$ from the recent work of Goncharov and Rudenko that was used in their proof of the weight 4 case of Zagier's Polylogarithm Conjecture.
	\end{abstract}
\maketitle
	
\section{Introduction}
In their recent breakthrough paper~\cite{goncharov2018motivic} Goncharov and Rudenko envisaged a very promising new strategy to prove Zagier's Polylogarithm Conjecture (ZPC) by relating it to cluster algebra complexes. A crucial ingredient in their proof of the case of the conjecture in weight 4 was a new functional equation relating the multiple polylogarithms $I_{3,1}(x,y)$ and $\Li_4(z)$, which they denoted by ${\Qb_4}$.

This result prompted our experimental search for higher analogues, and our computer implementation allowed us to find analogues, with a combinatorial structure inspired by and quite reminiscent of their ${\Qb_4}$, for higher weights.

These findings date back to 2018 at the MPI Bonn and were both communicated to Don Zagier and subsequently presented at a workshop on cluster algebras and the geometry of scattering amplitudes at the Higgs Centre in Edinburgh in March 2020. Our more ambitious goal of finding a bootstrapping procedure that would produce analogous results for general weight has now apparently been superseded by Rudenko's beautiful new preprint~\cite{rudenko2020depth} pertaining to Goncharov's depth conjecture.
As our approach does not seem to take the exact same symmetries into account, it has the big disadvantage of being harder to generalise but on the other hand it may have produced identities that are of a slightly different nature--both potentially in the use of symmetries and of the choice of functions---than the ones that we anticipate to appear eventually in his already announced `cluster polylogarithm' preprint (with Matveiakin).

We provide functional equations of type $\Qb_n$ up to weight $n=7$. Our formulas, found with the help of intensive computer calculations,  will likely differ from the ones Rudenko derives in that we use a potentially different  set of functions and impose cyclic symmetry. We anticipate that they could still be beneficial, in particular in view of the remaining task of `conditional' further depth reduction which may well result from combining  suitable specializations of the functional equations we give.

\section{Analogues of the functional equations $\Qb_n$ for  $n=5,6,7$} 
\label{sec:notation}
In \cite{Go1}, having sucessfully solved the weight~3 case of ZPC, Goncharov reduced the weight~4 case of his more encompassing Freeness Conjecture to an explicit calculation which would express $I_{3,1}(V(x,y),z)$, with the five term relation $V(x,y)$ in one slot, in terms of $\Li_4$. This was indeed shown to hold with 122 rather non-obvious terms (concocted as products of up to four cross ratios in 6 variables) in \cite{gangl2016multiple}. 
 
In \cite{goncharov2018motivic} Goncharov and Rudenko found an alternative and more conceptual way, introducing complexes of cluster algebras, to derive an equation that solves the same question without the need of giving those 122 terms explicitly,  and which furthermore has the important property of suggesting generalisations to higher weight.
The ensuing connection to moduli spaces ${\Mf}_{0,k}$ suggests to consider polyangulations of convex $2N$-gons for suitable integers $N$, resulting in pictures of the type given in the figures below. 

Further to the notation used in \cite{goncharov2018motivic} we introduce the following shorthand:
Any subpolygon comes equipped with a partition of its internal angles into two subsets (these often correspond in \cite{rudenko2020depth} to `even' and `odd' polytopes, but we note that our conventions allow successive even or odd indices for lower depth terms). We equip the angles of one of the two sets with a slice of pi(e), indicating that the associated argument is given as the cyclic ratio (already used extensively in~\cite{goncharov2018motivic})
	$$[x_{i_1},\dots,x_{i_{2m}}] = (-1)^m\frac{
	(x_{i_1}-x_{i_2})(x_{i_3}-x_{i_4})\dots (x_{i_{2m-1}} - x_{i_{2m}})} {(x_{i_2}-x_{i_3})(x_{i_4}-x_{i_5})\dots (x_{i_{2m}} - x_{i_{1}})}$$
where $i_1$ corresponds to any one of the $m$ labeled angles.
Moreover, we indicate the order in which the arguments are to be taken by indices 1, 2, etc. inside the given subpolygons. Finally, each picture stands for the sum over all cyclic permutations of the indices (corresponding to the rotations of the polygon).

Our functions are slight variants of the standard polylogarithm functions, indicated by an `$N$' in the notation, e.g. \( \IIr{3,1,1} \) stands for 
	$$ \IIr{3,1,1}(x,y,z) = -\Li_{3,1,1}((xyz)^{-1}, z,y)$$
and more generally
   \begin{align}
    \begin{split}
 \IIr{n_1\dots n_d}(a_1,\dots,a_d)
 &= I_{n_1,\dots,n_d}(a_1,(a_2\dots a_d)^{-1},(a_2\dots a_{d-1})^{-1},\dots,a_2^{-1}) \\
    &= (-1)^d \Li_{n_1,\dots,n_d}((a_1\dots a_d)^{-1},a_d,a_{d-1},\dots,a_2) 
    \end{split}
    \end{align}
(these variants have been chosen as they  satisfy simpler coproduct expressions than the underlying iterated integrals, see e.g. \cite[p.~248]{Ch}).
  
With the above explanations, we can now give the functional equations for weights 5, 6 and 7 that should play the role of the crucial  equation $\Qb_4$ in \cite{goncharov2018motivic}. 
The main characterising feature of these identities (and the reason why we propose them as generalisations of $\Qb_4$) is that the arguments of the iterated integrals are cross-ratios and higher cyclic-ratios attached to \( 2N \)-gons that form a complete polyangulation of a subpolygon by even polygons.  Moreover, the highest depth terms in these identities consist of a single depth \( d \) term of the form \( I_{d-n,\{1\}^n} \), evaluated at all possible full quadrangulations of a \( 2N \)-gon, up to cyclic symmetry.

\section{Identity in weight 5}
\label{sec:wt5}

\begin{figure}[h!]
\centering
\includegraphics[width=\textwidth]{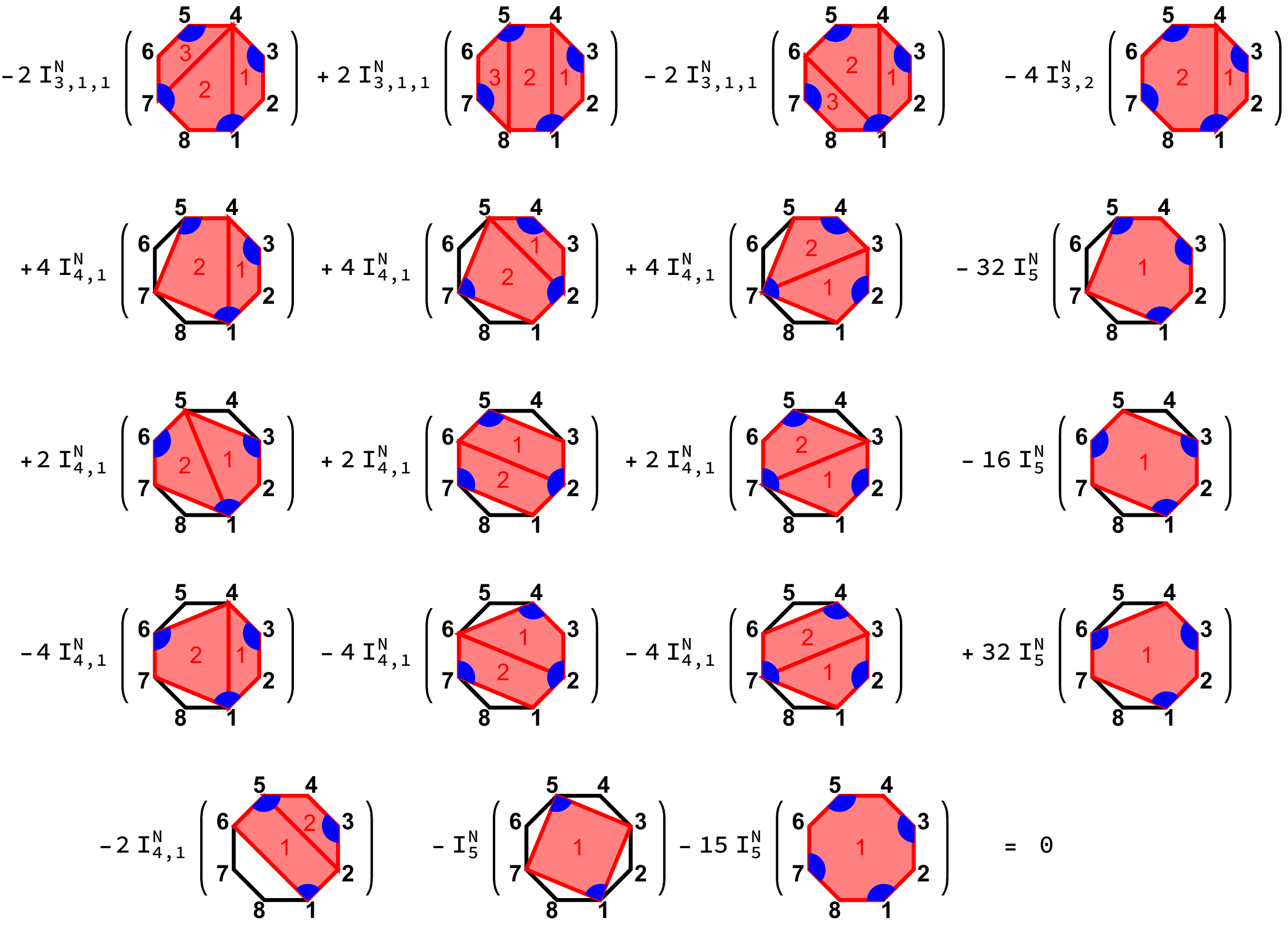}
\caption{A version of ${\bf Q}_5$}
\label{fig:q5}
\end{figure}

Our candidate for the functional equation~$\Qb_5$ is shown
in Figure~\ref{fig:q5}. It involves four different iterated integrals: 
$\IIr{3,1,1}$, $\IIr{3,2}$, $\IIr{4,1}$, and $\IIr{5}$.
\begin{Thm} \label{thm:wt5id}
	The weight $5$ multiple polylogarithm identity shown in 
	Figure~\ref{fig:q5} holds modulo products for generic $x_1,\dots,x_8$.
\end{Thm}

As explained in Section~\ref{sec:notation} each term in this figure
stands for cyclically symmetric sum of eight terms with indices taken modulo 8,
e.g.~the first term stands for the following expression
$$-4 \sum_{j=1}^8 \IIr{3,1,1}([x_{j+1},x_{j+2},x_{j+3},x_{j+4}], [x_{j+1}, x_{j+4}, x_{j+5}, x_{j+6}], [x_{j+1},x_{j+6}, x_{j+7},x_{j+8}]).$$

By a suitable specialization (or, more precisely, degeneration, since one needs to take limiting values) of the arguments we obtain an explicit reduction of a single term $I_{311}(x,y,z)$ to multiple polylogarithms of depth at most~2.

\begin{Prop}
	\label{prop:i311}
	The function $I_{3,1,1}(x,y,z)$ can be written as a linear combination of at most 47 terms, in terms of the functions \( I_{3,2} \) (1 term), \( I_{4,1} \) (23 terms), and \( I_5 \) (23 terms).

	\begin{proof}
		Specializing the identity of Theorem~\ref{thm:wt5id} to \( x_5=x_3=x_1 \) produces a combination of one generic \( \IIr{3,1,1} \) term, and several degenerate ones, up to the inversion and reversion of \( \IIr{3,1,1} \).  The generic term can then be isolated by subtracting the specialization of Theorem~\ref{thm:wt5id} to \( x_2=x_1 \).
	\end{proof}
\end{Prop}

\begin{Prop} \label{prop:i32}
	The function $I_{3,2}(x,y)$ can be written as a linear combination of at most 24 terms, in terms of the functions \( I_{4,1} \) (12 term), and \( I_5 \) (12 terms).
	
	
	\begin{proof}
		Specializing the above identity for \( \IIr{3,1,1} \) to \( x_7=x_2 \) reduces the \( \IIr{3,1,1} \) term to \( \IIr{3,2} \).  (The \( \IIr{3,2} \) term that occurs in the formula of Proposition~\ref{prop:i311} goes away under this degeneration.)
	\end{proof}
\end{Prop}

A different reduction of \( I_{3,2} \) to \( I_{4,1} \) and \( \Li_5\), involving a lot more complicated \( \Li_5 \) terms, but fewer and simpler \( I_{4,1} \) terms was given in \cite{Ch}.

\begin{Cor}
	Every weight 5 multiple polylogarithm can be expressed in terms of \( I_{4,1} \) and \( I_5 \).
	
	\begin{proof}
		This follows from the two preceding propositions together
		with an explicit reduction of any iterated integral in weight~$5$
		to $I_{3,1,1}$ that was given in~\cite{Ch}.
	\end{proof}
\end{Cor}

\section{Identities in weight 6}
\label{sec:wt6}
Our candidate for the identity~$\Qb_6$ is given
in Figure~\ref{fig:q6} below. It involves $\IIr{4,1,1}$, $\IIr{4,2}$, $\IIr{5,1}$, and $\IIr{6}$.
\begin{Thm} \label{thm:wt6id}
	The weight $6$ multiple polylogarithm identity shown in 
	Figure~\ref{fig:q6} holds modulo products for generic $x_1,\dots,x_{10}$.
\end{Thm}

\begin{figure}[h!]
	\centering
	\includegraphics[width=\textwidth]{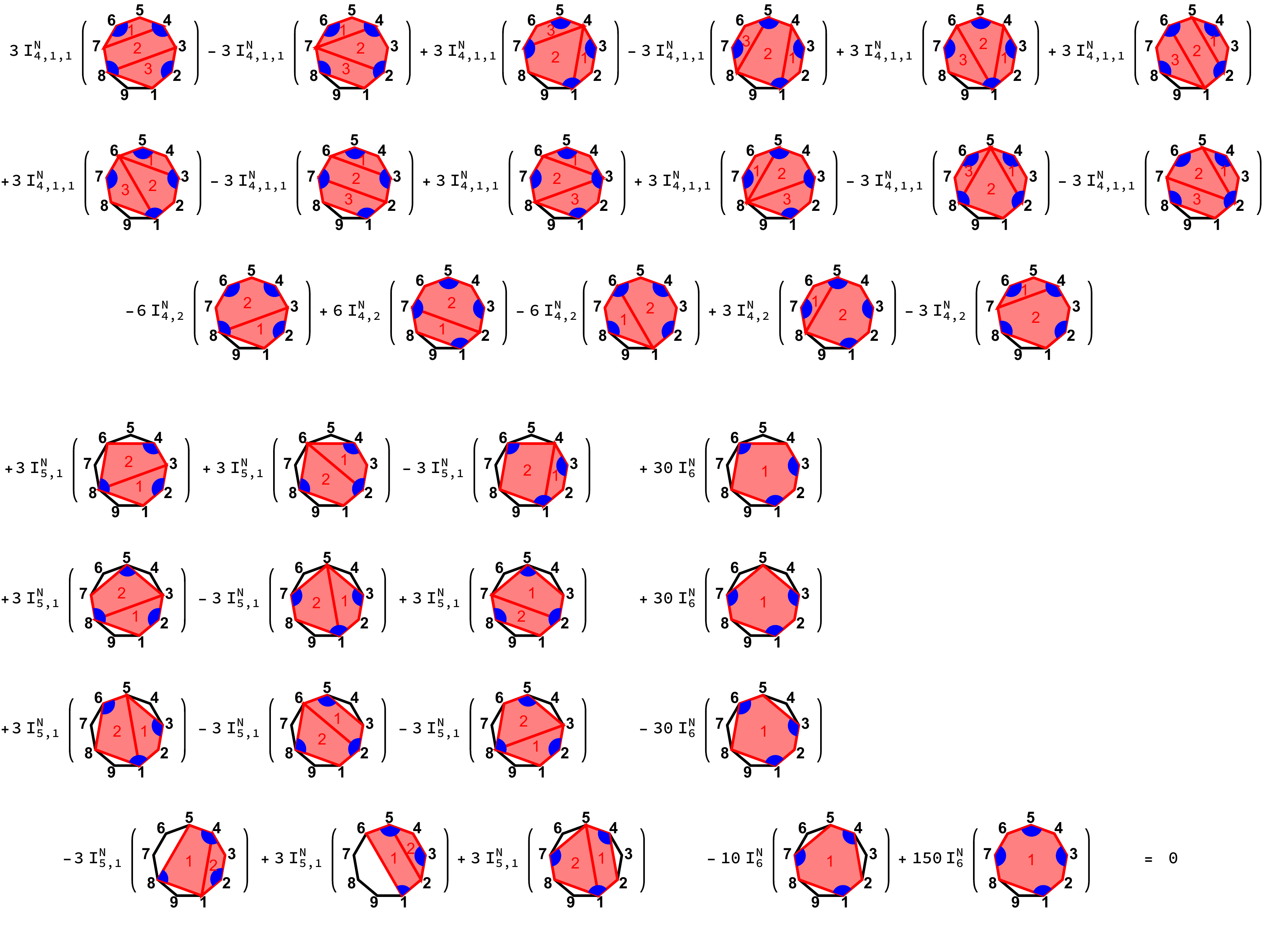}
	\vspace{-2.5em}
	\caption{A version of ${\bf Q}_6$}
	\label{fig:q6}
\end{figure}

\subsection{Reduction from depth~4 to depth~3 and identities of type $\Qb_n^{d}$}
To perform the depth reduction from 4 to 3 in the weight~6 case,
by analogy to what was done in the previous section, we need a slightly different, `off-diagonal' identity of type \( \Qb_n \).
The `diagonal' identities $\Qb_6$ and $\Qb_5$ presented above, as well as the 
original identities $\Qb_4$ and $\Qb_3$ of Goncharov and Rudenko have the property
that the highest depth of the iterated integrals that occur in $\Qb_n$
is $\lceil n/2 \rceil$. The weight $6$ identity whose top layer structure is indicated in Figure~\ref{fig:q6_4} below, however, has highest depth 4,
and we will call it $\Qb_{6}^{4}$, indicating the highest depth in the superscript. (In this notation the candidate for $\Qb_6$ given above would be
called $\Qb_{6}^{3}$.)

The specific identity that we give involves the functions
$\IIr{3,1,1,1}$, $\IIr{4,1,1}$, $\IIr{3,2,1}$, $\IIr{5,1}$, $\IIr{4,2}$, $\IIr{3,3}$, and $\IIr{6}$. We note that, unlike the identities
shown in Figure~\ref{fig:q5}, ~\ref{fig:q6}, and~\ref{fig:q7} (and breaking the convention that was set up in Section~\ref{sec:notation})  in this identity each term represents a \textit{signed} cyclic symmetrization, introducing a sign \( (-1)^j \) after cyclically shifting by \( j \) steps, i.e.
we replace each term by $\sum_{j=1}^{12}(-1)^jF(x_{1+j},\dots,x_{12+j})$ instead of
$\sum_{j=1}^{12}F(x_{1+j},\dots,x_{12+j})$.

\begin{Thm} \label{thm:wt6_4id}
	There exists a cyclically symmetric multiple polylogarithm 
	identity in weight~$6$ and depth~$4$ with~$168$ cyclic orbits of terms whose top layer structure is shown in Figure~\ref{fig:q6_4} that holds modulo products for generic $x_1,\dots,x_{12}$.
\end{Thm}

(We have also found a similar identity~$\Qb_{3}^{3}$ in weight 3 involving depth 3 functions and~$\Qb_{4}^{4}$ in weight~4 involving depth 4 functions
and also analogous identities~$\Qb_{5}^{3}$ and~$\Qb_{5}^{4}$ in weight~5.)

Using the identity~$\Qb_{6}^{4}$, similarly to the case of weight $5$ outlined in the previous section, one can obtain depth reduction of certain iterated integrals.

\begin{Prop}
	The function $I_{3,1,1,1}(x,y,z)$ can be written as a linear combination of the functions \( I_{3,2,1} \), \( I_{4,1,1} \), \( I_{4,2} \), \( I_{5,1} \), and  \( I_6 \).
	
	\begin{proof}
		Specializing to \( x_7=x_5=x_3=x_1 \) produces a combination consisting of one generic \( \IIr{3,1,1,1} \) term and several degenerate ones, up to the inversion and reversion relations for \( \IIr{3,1,1,1} \).  The generic term can then again be isolated by subtracting the specialization to \(x_4=x_1\).
	\end{proof}
\end{Prop}

\begin{Prop}\label{prop:i321}
	The function $I_{3,2,1}(x,y,z)$ can be expressed as a linear combination of \( I_{4,1,1} \), \( I_{4,2} \), \( I_{5,1} \) and \( I_6 \).
	
	\begin{proof}
		The proof goes along the same lines as the proof of Proposition~\ref{prop:i32}
		Specializing the above reduction for \( \IIr{3,1,1,1} \) to \(x_7=x_2\) reduces the \( \IIr{3,1,1,1} \) term to \( \IIr{3,2,1} \).  (The original \( \IIr{3,2,1} \) term goes away.) 
	\end{proof}
\end{Prop}
In general using the dihedral symmetries of \( I_{n_1,n_2,\ldots,n_k} \) allows one to express any iterated integral in terms $I_{3,\{1\}^a}$ or $I_{2,\{1\}^b,2,\{1\}^c}$.  By considering the stuffle product of \( \Li_{2,2,\{1\}^b} \) and \( \Li_{\{1\}^c} \) one can express \( I_{2,\{1\}^b,2,\{1\}^c} \) in terms of \( I_{2,\{1\}^{<b},2,\{1\}^{c'}} \) and lower depth.  Iteratively, this shows \( I_{2,2,\{1\}^{a-1}} \) suffices amongst integrals with indices involving \( 1 \)'s and two 2's.  Finally the shuffle product of \( I_{3,\{1\}^a} \) and \( I_{1} \) expresses \( I_{2,2,\{1\}^{a-1}} \) in terms of integrals with indices involving \( 1 \)'s and a single 3, meaning \( I_{3,\{1\}^a} \) alone is sufficient.
Combining this with Proposition~\ref{prop:i321} gives us.

\begin{Cor}
	Every weight 6 multiple polylogarithm can be expressed in terms of \( I_{4,1,1} \), \( I_{4,2} \), \( I_{5,1} \) and \( I_6 \).
\end{Cor}

\begin{figure}[h!]
	\centering
	\includegraphics[width=0.8\textwidth]{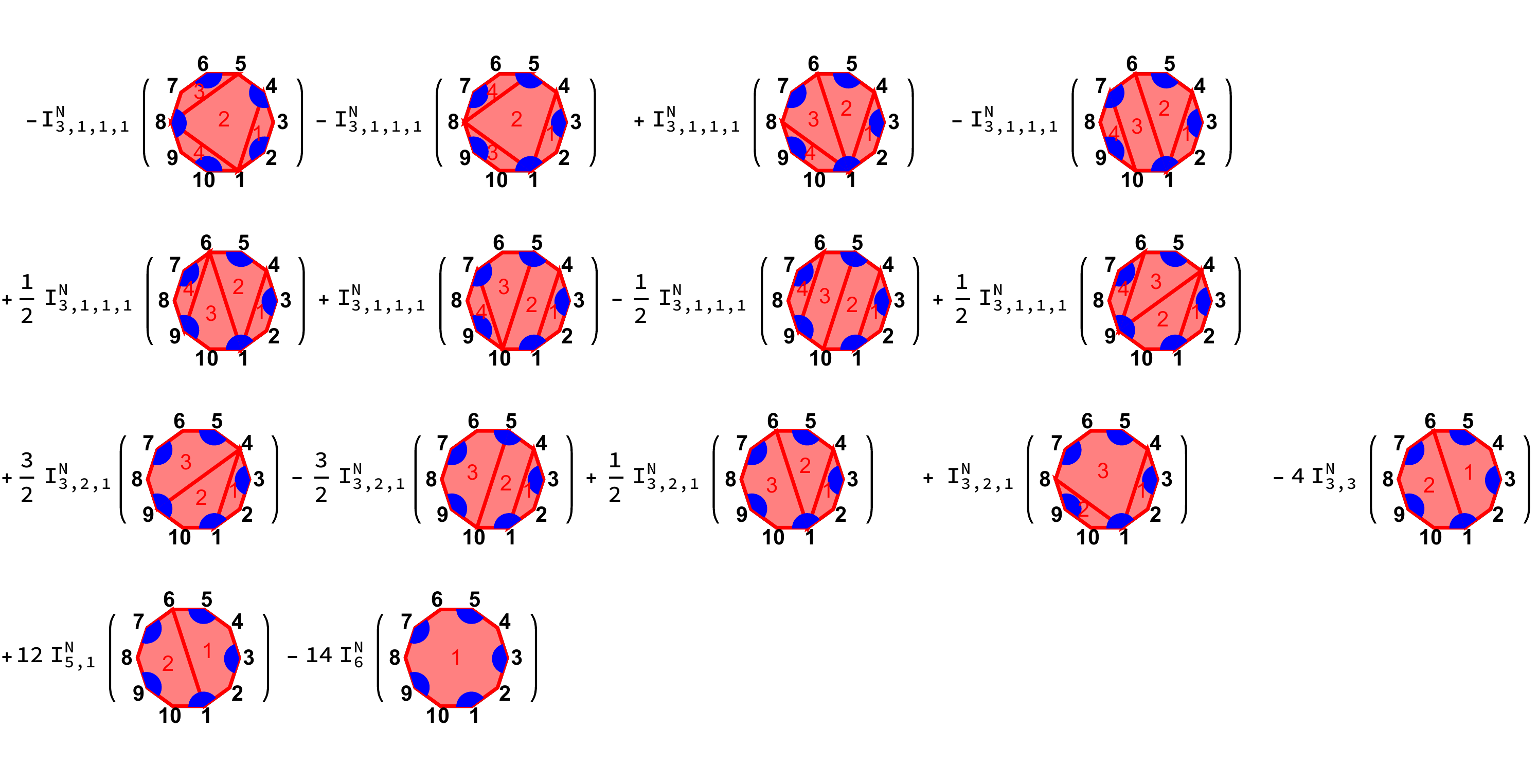}
	\vspace{-1.5em}
	\caption{Top layer of a weight 6 depth 4 identity}
	\label{fig:q6_4}
\end{figure}

\section{Identity in weight~7}
\begin{figure}[h]
	\centering
	\includegraphics[width=0.8\textwidth]{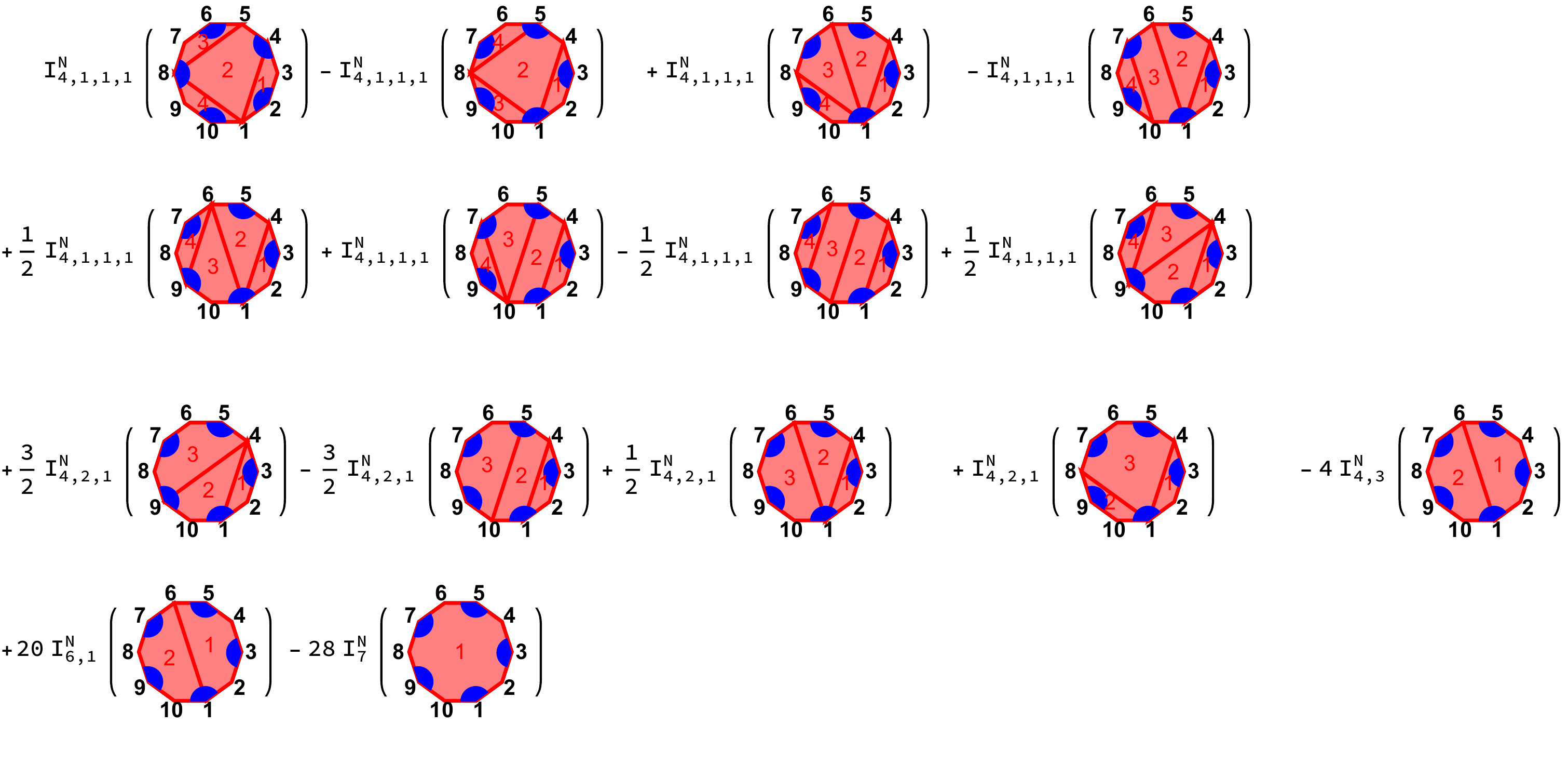}
	\vspace{-1.5em}
	\caption{Top layer of a version of $\Qb_7$}
	\label{fig:q7}
\end{figure}

Finally, we indicate the top level structure of our candidate for~$\Qb_7$ 
in Figure~\ref{fig:q7}. It involves the iterated integrals 
$\IIr{4,1,1,1}$, $\IIr{5,1,1}$, $\IIr{4,2,1}$, $\IIr{6,1}$, $\IIr{5,2}$, $\IIr{4,3}$, and $\IIr{7}$.

\begin{Thm} \label{thm:wt7id}
	There exists a cyclically symmetric multiple polylogarithm 
	identity in weight~$7$ and depth~$4$ with~$121$ cyclic orbits of terms whose top layer structure is shown in Figure~\ref{fig:q7} that holds modulo products for generic $x_1,\dots,x_{12}$.
\end{Thm}

Similarly to the proof of Proposition~\ref{prop:i311} we get the following.
\begin{Prop}
	The function $I_{4,1,1,1}(x,y,z)$ can be written as a linear combination of lower depth functions.
	\begin{proof}
		The same specialization to \( x_7=x_5=x_3=x_1 \) produces a combination of one generic \( \IIr{4,1,1,1} \) term, and several degenerate ones, up to the inversion and reversion relations for \( \IIr{4,1,1,1} \).  The generic term is isolated by subtracting the specialization to \( x_4=x_1 \).
	\end{proof}
\end{Prop}

\begin{Rem}
	It was checked with extensive linear algebra that \( I_{5,1,1} \) functions of the type appearing in \( \Qb_7 \), along with similar functions of lower depth, were sufficient to span the space of weight 7 iterated integrals.  More precisely, we were able to confirm that the rank of such functions agreed with the dimension computed by Brown in \cite{Br-rep}, for weight \( k = 7 \) iterated integrals in \( n = 8 \) projective variables (the dimension is \( \tfrac{1}{k} \sum_{d \mid k} \mu(\tfrac{k}{d}) \sum_{i=2}^{n-2} i^d = 53\,820 \)).  By the known reduction of weight 7 integrals to depth~5 (hence involving \( 5 + 3 = 8 \) projective variables), this implies that every weight 7 integral can be expressed in terms of \( I_{5,1,1} \), \( I_{5,2} \), \( I_{6,1} \) and \( I_7 \).  Given the large dimension, we were unable to extract explicit formulas for this reduction.
\end{Rem}

\end{document}